# The exponential rank of nonarchimedean exponential fields

*Franz-Viktor and Salma Kuhlmann*[*]

16. 6. 1997

## 1 Introduction

Based on the work of Hahn, Baer, Ostrowski, Krull, Kaplansky and the Artin-Schreier theory, and stimulated by the paper [L] of S. Lang in 1953, the theory of real places and convex valuations has witnessed a remarkable development and has become a basic tool in the theory of ordered fields and real algebraic geometry. Surveys on this development can be found in [LAM] and [PC]. In this paper, we take a further step by adding an exponential function to the ordered field. Beforehand, let us sketch the basic facts about convex valuations.

Let $(K, <)$ be an ordered field and $w$ a valuation of $K$, with valuation ring $R_w$, valuation ideal $I_w$, value group $wK$ and residue field $Kw$. Then $w$ is called **compatible with the order** if and only if it satisfies, for all $x, y \in K$:

**(CO)**     $a \leq b \leq 0 \ \vee \ a \geq b \geq 0 \implies wa \leq wb$ .

We will denote the valuation ring of $w$ by $R_w$ and the valuation ideal by $I_w$. It is well known that (CO) is equivalent to each of the following assertions:
1) $R_w$ is a convex subset of $(K, <)$,
2) $I_w$ (or equivalently, the set $1 + I_w$ of 1-units) is a convex subset of $(K, <)$,
3) the positive cone of $(K, <)$ contains $1 + I_w$,
4) $I_w < 1$,
5) the image of the positive cone of $(K, <)$ under the residue map $K \ni a \mapsto aw \in Kw$ is a positive cone in $Kw$.

The last assertion means that the order of $K$ canonically induces an order on the residue field $Kw$. In view of 1) and 2), a valuation $w$ compatible with the order is also said to be a **convex valuation**. For every convex valuation $w$, the set $\mathcal{U}_w^{>0} := \{a \in K \mid wa = 0 \wedge a > 0\}$ of **positive units** of $R_w$ is a convex subgroup of the ordered multiplicative group $(K^{>0}, \cdot, 1, <)$ of positive elements of $K$.

[*]Supported by a Deutsche Forschungsgemeinschaft fellowship

Assume that $w$ is a **coarsening** of a second valuation $v$, that is, $R_w \supset R_v$ (in this case, we also say that $v$ is **finer** than $w$). Then $I_w \subset I_v$, and if $v$ is convex, it follows from condition 3) or 4) that also $w$ is convex. There is always a finest convex valuation $v$, and the other convex valuations are precisely the coarsenings $w$ of $v$. This valuation $v$ is called the **natural valuation of** $(K, <)$.

Let us quickly recall how natural valuations are obtained already on a totally ordered abelian group $(G, <)$. We set $|a| := \max\{a, -a\}$. Two elements $a, b$ are **archimedean equivalent** if there is some $n \in \mathbb{N}$ such that $n|a| \geq |b|$ and $n|b| \geq |a|$. The equivalence class of $a$ is called **archimedean class of** $a$ and is denoted by $[a]$. The set $\{[a] \mid 0 \neq a \in G\}$ is totally ordered by setting $[a] < [b]$ if and only if $|a| > n|b|$ for all $n \in \mathbb{N}$. The map $v : a \mapsto [a]$ is the **natural valuation** of $G$. It satisfies the triangle inequality $v(a+b) \geq \min\{va, vb\}$ and $v(-a) = va$ as well as property (CO). The natural valuation of an ordered field $(K, <)$ is just the natural valuation on its ordered additive group $(K, +, 0, <)$. In this case, $vK := v(K \setminus \{0\})$ forms a totally ordered abelian group endowed with the addition $[a] + [b] := [ab]$, and $v$ is a field valuation. It is characterized by the fact that its residue field is an archimedean ordered field. The valuation ideal $I_v$ is the set of all **infinitesimals**, and $K^{>0} \setminus R_v$ is the set of all **positive infinite elements** (elements which are incomparably bigger than 1). For more information on natural valuations, see [KS].

Throughout this paper, $K$ will be a nonarchimedean ordered field, and $v$ will denote its nontrivial natural valuation. The natural valuations of the appearing ordered groups will be denoted by $v_G$, and $v_G G$ shall be the ordered set $v_G(G \setminus \{0\})$.

By general valuation theory, the set $\mathcal{R}$ of all valuation rings $R_w$ of coarsenings $w \neq v$ of $v$ is totally ordered by inclusion, and it is order isomorphic to the set of all nontrivial convex subgroups of the value group $vK$ (again ordered by inclusion). Its order type is called the **rank** of $v$ (or in our case, of $(K, <)$); for convenience, we will identify it with $\mathcal{R}$. The convex subgroup corresponding to $R_w$ is

$$G_w := \{va \mid a \in K \wedge wa = 0\} = v(\mathcal{U}_w^{>0}),$$

and the value group $wK$ is canonically isomorphic to $vK/G_w$. For example, the rank of an archimedean ordered field is empty since its natural valuation is trivial. The rank of the rational function field $K = \mathbb{R}(t)$ with any order is a singleton, $\mathcal{R} = \{K\}$, and $vK$ is the only nontrivial convex subgroup of the nontrivial archimedean value group $vK$. The reader may have noted that we are not using the classical definition of "rank" since we include the trivial valuation but exclude $v$. The present version will be more useful for our purposes.

Now we add an exponential to the ordered field. For an **exponential** $f$ on $(K, <)$ we only require that it is an isomorphism from the ordered additive group $(K, +, 0, <)$ onto $(K^{>0}, \cdot, 1, <)$. Let $w$ be a convex valuation on $K$. Then we will say that $w$ and $f$ are **compatible** if the following holds:

**(CE)** $\qquad f(R_w) = \mathcal{U}_w^{>0}$ and $f(I_w) = 1 + I_w$.

Since $Kw = R_w/I_w$ and $(Kw^{>0}, \cdot, 1, <) = (\mathcal{U}_w^{>0}, \cdot, 1, <)/1 + I_w$, this means that $f$ induces canonically an exponential $fw : (Kw, +, 0, <) \to (Kw^{>0}, \cdot, 1, <)$ on the residue field $Kw$.

(This is the analogue to the characteristic property of convex valuations to canonically induce an order on the residue field.)

Let us mention that if $K$ is a model of the elementary theory $T$ of an exponentially bounded o-minimal expansion of the reals, such that the exponential $f$ is definable, then the valuation rings $R_w$ of valuations $w$ compatible with $f$ are precisely the $T$-convex valuation rings of $K$, in the sense of [DL].

The valuation rings $R_w$ of nontrivial convex valuations $w$ satisfying the first condition of (CE) form a subset $\mathcal{R}_f$ of $\mathcal{R}$. Its order type will be called the **exponential rank** of the exponential field $(K, <, f)$; again, we identify it with $\mathcal{R}_f$. We wish to characterize the corresponding convex subgroups $G_w$. However, these subgroups do not carry any information concerning the second condition of (CE). So it may well happen that $\mathcal{R}_f$ also contains valuation rings of valuations which are not compatible with $f$. A natural way to overcome this deficiency is to require that $f$ satisfies the following elementary axiom:

**(T$_1$)** $\qquad v(f(a) - 1 - a) > va \qquad$ for all $a \in I_v$.

It belongs to a scheme of axioms which gives a valuation theoretical interpretation of the Taylor expansion of the usual exponential function on $\mathbb{R}$ (see [KK1] for details). If the natural valuation $v$ is compatible with $f$ and (T$_1$) holds, then we call $f$ a **T$_1$-exponential**. Note that if an ordered field $K$ admits any exponential, then it admits an exponential compatible with the natural valuation (cf. [KS], Section 3.3). If $f$ is a T$_1$-exponential, then $f(I_w) = 1 + I_w$ holds for every coarsening $w$ of $v$; this is a consequence of Lemma 7 in the next section. Then $f$ is compatible with $w$ if and only if $f(R_w) = \mathcal{U}_w^{>0}$ and consequently, $\mathcal{R}_f$ is precisely the set of all valuation rings of valuations which are compatible with $f$.

We shall characterize the subgroups $G_w$ of $vK$ for which $R_w \in \mathcal{R}_f$ by use of a **contraction map** $\chi$ induced on $vK$ by the exponential $f$ (more precisely, by its inverse, the **logarithm** $\ell$). This map in turn induces a map $\zeta$ on the rank $\mathcal{R}$. For the details, see Section 3. To avoid unpleasant case distinctions which would make the theory complicated without telling anything more about the interesting cases, we fix the "orientation" of these two maps. This is done by requiring that $f$ satisfies the following elementary **growth axiom** scheme, which is also satisfied by the usual exponential function on $\mathbb{R}$:

**(GA)** $\qquad a > n^2 \implies f(a) > a^n \qquad (n \in \mathbb{N})$.

If this holds, then $f$ will be called a **strong exponential**. The following gives a basic characterization of the convex valuations which are compatible with such an exponential:

**Theorem 1** *Let $f$ be a strong $T_1$-exponential. A coarsening $w$ of $v$ is compatible with $f$ if and only if for every $a \in K$,*

$$va \in G_w \Rightarrow vf(a) \in G_w. \tag{1}$$

The proof and further characterizations by use of the maps $\chi$ and $\zeta$ will be given in Section 3. There, we will also introduce equivalence relations $\sim_\chi$ and $\sim_\zeta$ induced by $\chi$ and $\zeta$ (in the spirit of archimedean equivalence).

The convex subgroup $G_w$ of $G$ is called **principal** if there is some $g \in G$ such that $G_w$ is the minimal convex subgroup containing $g$ (it exists since the intersection of all convex subgroups containing $g$ is a convex subgroup). We will call $G_w$ $f$-**principal** if there is some $g \in G$ such that $G_w$ is the minimal convex subgroup containing $g$ and closed under the map $va \mapsto vf(a)$ (this definition works for arbitrary maps $f$ on $K$). By the **principal rank** of $(K, <)$ we mean the subset $\mathcal{R}^{\mathrm{pr}}$ of $\mathcal{R}$ consisting of all $R_w \in \mathcal{R}$ for which $G_w$ is principal. If $f$ is a strong exponential, then the **principal exponential rank** shall be the subset $\mathcal{R}^{\mathrm{pr}}_f$ of $\mathcal{R}_f$ consisting of all $R_w \in \mathcal{R}_f$ for which $G_w$ is $f$-principal. Note that by the above theorem, a strong $T_1$-exponential is compatible with $w$ already if $G_w$ is $f$-principal. By induction, we define $f^1(a) := f(a)$ and $f^{n+1}(a) := f(f^n(a))$.

**Theorem 2** *Let $R_w \neq R_v$ be a convex valuation ring of the nonarchimedean ordered field $(K, <)$.*
*a) $R_w$ lies in the principal rank if and only if for some $b$ in the residue field $Kw$, the sequence $(b^n)_{n \in \mathbb{N}}$ is cofinal in $Kw$. The principal rank is a singleton (namely, $\mathcal{R}^{\mathrm{pr}} = \{K\}$) if and only if $(a^n)_{n \in \mathbb{N}}$ is cofinal in $K$ for every positive infinite element $a \in K$.*
*b) Assume in addition that $f$ is a strong $T_1$-exponential and that $w$ is compatible with $f$. Then $R_w$ lies in the principal exponential rank if and only if for some $b$ in the residue field $Kw$, the sequence $((fw)^n(b))_{n \in \mathbb{N}}$ is cofinal in $Kw$. The principal exponential rank is a singleton (namely, $\mathcal{R}^{\mathrm{pr}}_f = \{K\}$) if and only if $(f^n(a))_{n \in \mathbb{N}}$ is cofinal in $K$ for every positive infinite element $a \in K$.*

Note that the valuation ring $R_v$ of the natural valuation $v$ would lie in the principal rank as well as in the principal exponential rank if we would not exclude it from the rank $\mathcal{R}$. In contrast to this, $K$ lies in the principal rank if and only if there is some $a \in K$ such that $(a^n)_{n \in \mathbb{N}}$ is cofinal in $K$, and it lies in the principal exponential rank if and only if there is some $a \in K$ such that $(f^n(a))_{n \in \mathbb{N}}$ is cofinal in $K$. Therefore, we want to include $K$ in the rank. But then, we have to exclude $R_v$ since an ordered field $(K, <)$ having precisely $n$ convex valuation rings $\neq K$ is said to have rank $n$, which should correspond to the number of elements in $\mathcal{R}$.

The description of the set of all liftings of an order through a place is well known. In the same spirit, we will describe in Theorem 16 the set of all liftings of a logarithm. We work with logarithms rather than with exponentials since we can give this description even for **non-surjective logarithms**, i.e., embeddings of $(K^{>0}, \cdot, 1, <)$ in $(K, +, 0, <)$. These will play a crucial role in the following.

So far, we have only described results that are in nice analogy to the theory of real places. But when it comes to existence results, the analogy breaks down. If a field has a place onto an ordered residue field, then the order can be lifted up to the field through the place. It is not surprising that exponentials cannot be lifted through arbitrary places. But one might expect that certain closure properties (like "henselian place", "divisible value group" or perhaps some compatibility of the value group) would make such a lifting possible. For example, if $k$ is an ordered field and $G$ a nontrivial ordered abelian group, then the ("generalized") power series field $K = k((G))$ admits at least one nonarchimedean order. Further, $K$ is real closed if and only if $k$ is real closed and $G$ is divisible. This provides a simple and elegant method of constructing nonarchimedean ordered real closed fields of any given rank $\tau$, provided that we can construct a divisible ordered abelian group of rank $\tau$. But the latter is easy: we just take $G$ to be the lexicographic sum of copies of $\mathbb{Q}$ with index set $\tau$ (or the corresponding Hahn product if we want to stick to the philosophy of power series). If $k$ is archimedean, the principal rank of $K$ will then be $\tau$. But for the construction of exponential fields with arbitrary given principal exponential rank, this approach fails. Indeed, we have shown in [KKS] that power series fields *never* admit exponentials compatible with their canonical valuation (and if $k$ is archimedean, then there is no exponential on $k((G))$ at all).

Nevertheless, we are able to construct exponential fields with arbitrary principal exponential rank. This is done in two steps. First, we construct non-surjective logarithms on power series fields. Therefore, we have to develop our theory of exponential rank and define the maps $\chi$ and $\zeta$ also for such logarithms, provided they satisfy adequate versions of the axioms (T$_1$) and (GA). This is done in Section 3. In the second step, we obtain a surjective logarithm by taking the union over a suitable countable ascending chain of such power series fields with non-surjective logarithms. This is done in Section 6. For the case of models of the theory $T_{\mathrm{an}}(\exp)$ of the reals with restricted analytic functions and exponential function (cf. [DMM1]), we shall prove:

**Theorem 3** *Take a model of $T_{\mathrm{an}}(\exp)$ and an order type $\tau$ which extends the principal exponential rank of $(K, f)$. Then $(K, f)$ can be elementarily embedded in a model $(K_\omega, f_\omega)$ of $T_{\mathrm{an}}(\exp)$ which is a countable union of power series fields and has principal exponential rank $\tau$. The embedding can be chosen to be truncation closed.*

"Truncation closed" means that the truncation of any power series in the image of the embedding lies again in this image. Note that as we exclude $R_v$ from the principal exponential rank, we do not have to require that $\tau$ has a smallest element.

Our construction given in Section 6 for the proof of Theorem 3 is rather abstract. In [KK3], we give an explicit construction, which helps to exhibit the connection between endomorphisms of the exponential rank and the growth rate of the constructed exponentials. This allows us to construct on a fixed real closed field infinitely many exponentials of distinct exponential rank. Thus, in contrast to the rank, the exponential rank of a real closed exponential field is in general not uniquely determined.

## 2 Strong logarithms and $T_1$-logarithms

If the logarithm $\ell$ is surjective, then the exponential $\ell^{-1}$ is compatible with $w$ if and only if $\ell(\mathcal{U}_w^{>0}) = R_w$ and $\ell(1 + I_w) = I_w$. Thus, for a not necessarily surjective logarithm $\ell$ we will say that $w$ **and $\ell$ are compatible** if

$$\ell(\mathcal{U}_w^{>0}) = R_w \cap \mathrm{im}(\ell) \quad \text{and} \quad \ell(1 + I_w) = I_w \cap \mathrm{im}(\ell) . \tag{2}$$

We let $\mathcal{R}_\ell$ denote the subset of $\mathcal{R}$ containing all $R_w$ for which $w$ satisfies the first condition of (2), and we call it the **exponential rank of $\ell$**. If $\ell = f^{-1}$, then $\mathcal{R}_\ell = \mathcal{R}_f$. Note that the first condition of (2) is equivalent to:

$$\ell(\mathcal{U}_w^{>0}) \subset R_w \quad \text{and} \quad \ell(K^{>0} \setminus \mathcal{U}_w^{>0}) \subset K \setminus R_w . \tag{3}$$

We will now consider the content of the axioms (GA) and (T$_1$) for logarithms. Because of the condition "$a > n^2$", axiom scheme (GA) is void for infinitesimals. That is, it gives information only in the case of $va \leq 0$. It holds for $va = 0$ if the exponential $fv$ induced by $f$ on $Kv$ satisfies (GA) (e.g. if $Kv = \mathbb{R}$ and $fv$ is the usual exponential function); the proof is simple, see e.g. [KK1], Lemma 2.10.

Now we are interested in the case of $va < 0$. In this case, "$a > n^2$" holds for all $n \in \mathbb{N}$ if $a$ is positive. Restricted to $K \setminus R_v$, axiom scheme (GA) is thus equivalent to the assertion

$$\forall n \in \mathbb{N}: \ f(a) > a^n \qquad \text{for all } a \in K^{>0} \setminus R_v . \tag{4}$$

Applying the logarithm $\ell = f^{-1}$ on both sides, we find that this is equivalent to

$$\forall n \in \mathbb{N}: \ a > \ell(a^n) = n\ell a \qquad \text{for all } a \in K^{>0} \setminus R_v . \tag{5}$$

Via the natural valuation $v$, this in turn is equivalent to

$$va < v\ell a \qquad \text{for all } a \in K^{>0} \setminus R_v . \tag{6}$$

A logarithm $\ell$ (whether surjective or not) will be called a **strong logarithm** if it satisfies this condition. A real closed ordered field $(K, <)$ with exponential $f$ is a model of real exponentiation if it is a model of restricted real exponentiation and $f$ satisfies (GA); this is the content of Ressayre's Theorem (cf. [RE]), and it also holds if one adds restricted analytic functions (cf. [DMM1], (4.10)). So let us note:

**Lemma 4** *Let $K$ be a model of the reals with restricted analytic functions, and $f$ an exponential on $K$ such that $f$ coincides on $[-1, 1]$ with the interpretation of the (symbol for the) restricted exponential. Then $(K, f)$ is a model of $T_{\text{an}}(\exp)$ if and only if $\ell = f^{-1}$ is a strong logarithm.*

Assume that $w$ is a coarsening of $v$. Then because of $R_w \supset R_v$ we have that $va < vb$ implies $wa \leq wb$. Hence, (6) implies:

$$wa \leq w\ell a \qquad \text{for all } a \in K^{>0} \setminus R_v . \tag{7}$$

This in turn implies that

$$\ell(\mathcal{U}_w^{>0} \setminus R_v) \subset R_w . \tag{8}$$

We note that

$$\mathcal{U}_w^{>0} \ = \ (\mathcal{U}_w^{>0} \cap I_v) \cup \mathcal{U}_v^{>0} \cup (\mathcal{U}_w^{>0} \setminus R_v) \quad \text{with} \quad (\mathcal{U}_w^{>0} \cap I_v) < \mathcal{U}_v^{>0} < (\mathcal{U}_w^{>0} \setminus R_v)$$

and $\mathcal{U}_w^{>0} \cap I_v = \{x^{-1} \mid x \in \mathcal{U}_w^{>0} \setminus R_v\}$. If $\ell$ is compatible with $v$, then $\ell(\mathcal{U}_v^{>0}) \subset R_v \subset R_w$. If $\ell$ also satisfies (8), then $\ell(\mathcal{U}_w^{>0} \cap I_v) = -\ell(\mathcal{U}_w^{>0} \setminus R_v) \subset -R_w = R_w$, and

$$\ell(\mathcal{U}_w^{>0}) \ = \ \ell(\mathcal{U}_w^{>0} \cap I_v) \cup \ell(\mathcal{U}_v^{>0}) \cup \ell(\mathcal{U}_w^{>0} \setminus R_v) \ \subset \ R_w . \tag{9}$$

Using this fact, we prove:

**Lemma 5** *Assume that $\ell$ is a (not necessarily surjective) strong logarithm which is compatible with $v$. Then the first condition of (2) is equivalent to:*

$$\ell(K^{>0} \setminus R_w) \subset K^{>0} \setminus R_w . \tag{10}$$

Proof: As the first condition of (2) is equivalent to (3) and as (9), the first condition of (3), holds under the assumptions of the lemma, it remains to show that the second condition of (3), $\ell(K^{>0} \setminus \mathcal{U}_w^{>0}) \subset K \setminus R_w$, is equivalent to (10). Since $wa > 0 \Leftrightarrow wa^{-1} < 0$ and $w\ell a^{-1} = w(-\ell a) = w\ell a$, it is equivalent to

$$a \in K^{>0} \wedge wa < 0 \Rightarrow w\ell a < 0 \,. \tag{11}$$

Note that $a \in K^{>0} \setminus R_w$ implies that $a > 1$ and thus, $\ell a > 0$. Therefore, (11) is equivalent to (10). □

We turn to the Taylor axiom $(T_1)$.

**Lemma 6** *If $\ell$ is the inverse of an exponential $f$ which is compatible with $v$, then the Taylor axiom $(T_1)$ is equivalent to*

$$v(b - \ell(1+b)) > vb \qquad \text{for all } b \in I_v \,. \tag{12}$$

Proof: By the compatibility, every $a \in I_v$ is of the form $\ell(1+b)$ with $b \in I_v$, and every such $\ell(1+b)$ is in $I_v$. With $a = \ell(1+b)$, the assertion $v(f(a) - 1 - a) > va$ is equivalent to $v(b - \ell(1+b)) > v\ell(1+b)$. But as this implies that $vb = v\ell(1+b)$, it is equivalent to $v(b - \ell(1+b)) > vb$. □

This leads to the following definition: a logarithm $\ell$ (whether surjective or not) will be called a **$(T_1)$-logarithm** if it is compatible with $v$ and satisfies condition (12).

**Lemma 7** *For a $(T_1)$-logarithm $\ell$, the condition $\ell(1 + I_w) = I_w \cap \mathrm{im}(\ell)$ holds for all coarsenings $w$ of $v$.*

Proof: Condition (12) implies that $vy = v\ell(1+y)$ for all $y \in I_v$ and therefore, that $\ell(1 + I_w) \subset I_w$ and $\ell(1 + I_v \setminus 1 + I_w) \subset I_v \setminus I_w$ for every coarsening $w$ of $v$. By definition, $\ell$ is compatible with $v$, so we have that $\ell(1 + I_v) = I_v \cap \mathrm{im}(\ell) \supset I_w \cap \mathrm{im}(\ell)$. Consequently, $\ell(1 + I_w) = I_w \cap \mathrm{im}(\ell)$. □

By this lemma, a $T_1$-logarithm $\ell$ always satisfies the second condition of (2). So we have proved:

**Lemma 8** *Let $\ell$ be a strong $T_1$-logarithm. Then a coarsening $w$ of $v$ is compatible with $\ell$ if and only if it satisfies (10). Further, the exponential rank $\mathcal{R}_\ell$ is the subset of all $R_w \in \mathcal{R}$ for which $w$ is compatible with $\ell$.*

## 3 The maps $\chi$ and $\zeta$

Throughout this section, we assume $\ell$ to be a strong logarithm, not necessarily surjective or $T_1$. We shall now study the structure which $\ell$ induces on the value group and the rank, and deduce conditions for a coarsening $w$ of $v$ to satisfy (10). We set $(wK)^{<0} := \{g \in wK \mid g < 0\}$ and note that $v_G wK = \{v_G g \mid 0 \neq g \in wK\} = \{v_G g \mid g \in (wK)^{<0}\}$.

• **Definition of the maps $\chi_w$ and $\zeta_w$ on the value group and rank.**

Suppose that $w$ satisfies (10). We define a map $\chi_w : (wK)^{<0} \to (wK)^{<0}$ by setting

$$\chi_w(wa) := w\ell a \qquad \text{for all } a \in K^{>0} \setminus R_w .$$

This definition does not depend upon the representative $a$ of the value $wa$. Indeed, if $a, b \in K^{>0} \setminus R_w$ such that $wa = wb$, then $a = bc$ with $c \in \mathcal{U}_w^{>0}$. It follows that $\ell a = \ell(bc) = \ell b + \ell c$, with $\ell c \in R_w$. Since $w\ell b < 0$, we obtain that $w\ell a = w(\ell b + \ell c) = \min\{w\ell b, w\ell c\} = w\ell b$.

Now suppose that $g = wa$ and $g' = wa'$ are elements of $(wK)^{<0}$ with $a, a' \in K^{>0} \setminus R_w$ and $0 > g > g'$. Then $g$ and $g'$ are archimedean equivalent if and only if there is $n \in \mathbb{N}$ such that $ng < g'$, that is, $wa^n < wa'$. In this case, $a^n > a' > a$. This implies that $n\ell a = \ell a^n > \ell a' > \ell a$, hence $w\ell a' = w\ell a$, showing that $\chi_w g = \chi_w g'$. That is, every two archimedean equivalent elements of $(wK)^{<0}$ have the same image under $\chi_w$.

As a consequence, the map $\zeta_w : v_G wK \to v_G wK$ given by

$$\zeta_w(v_G g) := v_G \chi_w g \qquad \text{for all } g \in (wK)^{<0}$$

is well defined. We set $\chi := \chi_v$ and $\zeta := \zeta_v$.

Since $\chi_w$ contracts the negative part of every archimedean class of $wK$ to one element, we call it a **contraction**. The properties and the model theory of such maps have been studied in detail in [KK1], [KF1] and [KF2].

By definition of $\chi_w$ and $\zeta_w$, we have the following commutative diagram:

$$\begin{array}{ccc}
K^{>0} \setminus R_w & \xrightarrow{\ell} & K^{>0} \setminus R_w \\
\downarrow w & /// & \downarrow w \\
(wK)^{<0} & \xrightarrow{\chi_w} & (wK)^{<0} \\
\downarrow v_G & /// & \downarrow v_G \\
v_G wK & \xrightarrow{\zeta_w} & v_G wK
\end{array}$$

In this diagram, the map $w$ reverses the order $\leq$, and $v_G$ preserves the order $\leq$. Moreover, both are onto. Since also $\ell$ preserves the order $\leq$ (i.e., $\ell$ is monotone), we find that:

a) $\chi_w$ and $\zeta_w$ are monotone,

b) if $\ell$ is onto, then so are $\chi_w$ and $\zeta_w$.

Since $\chi_w$ and $\zeta_w$ are in general not injective, they may not be strictly monotone.

(If $w$ is the trivial valuation, then $K^{>0} \setminus R_w$, $(wK)^{<0}$ and $v_G wK$ are empty and $\chi_w$ and $\zeta_w$ are the empty maps.)

Since $\chi_w$ preserves $\leq$ and sends archimedean equivalent elements (i.e., elements with equal $v_G$-value) to one point, the following holds:

$$\left. \begin{array}{l} v_G g = v_G g' \Rightarrow \chi_w g = \chi_w g' \\ v_G g \geq v_G g' \Rightarrow \chi_w g \geq \chi_w g' \end{array} \right\} \qquad \text{for all } g, g' \in (wK)^{<0} . \qquad (13)$$

From (6) and (7) we infer:

$$g < \chi g \text{ for all } g \in (vK)^{<0} \quad \text{and} \quad g \leq \chi_w g \text{ for all } g \in (wK)^{<0} . \qquad (14)$$

It follows that $v_G g \leq v_G \chi_w g$ for all $g \in (wK)^{<0}$. But $v_G g = v_G \chi g$ cannot hold; otherwise (13) would yield that $\chi g = \chi\chi g$, in contradiction to (14). So we find:

$$\gamma < \zeta\gamma \quad \text{and} \quad \gamma \leq \zeta_w \gamma \qquad \text{for all } \gamma \in v_G wK \, . \tag{15}$$

- **Equivalence relations induced by $\ell$, $\chi_w$ and $\zeta_w$.**

If $\varphi$ is any map from a totally ordered set $S$ into itself, we define a relation $\sim_\varphi$ on $S$ by setting $a \sim_\varphi a'$ if the convex hulls of $\{a, \varphi^n(a) \mid n \in \mathbb{N}\}$ and $\{a', \varphi^n(a') \mid n \in \mathbb{N}\}$ have a nonempty intersection. This relation is in general not transitive. But if $\varphi$ is monotone, it is an equivalence relation. In this case, we will say that $a$ **and $a'$ are $\varphi$-equivalent** if $a \sim_\varphi a'$. The equivalence classes $[a]_\varphi$ of $\sim_\varphi$ are convex and closed under application of $\varphi$. By the convexity, the order of $S$ induces an order on $S/\sim_\varphi$ such that $[a]_\varphi < [b]_\varphi$ if and only if $a' < b'$ for all $a' \in [a]_\varphi$ and $b' \in [b]_\varphi$. On the positive part or the negative part of an ordered abelian group $G$, the archimedean equivalence relation is obtained by setting $\varphi(a) := 2a$, and $v_G$ is the map $a \mapsto [a]_\varphi$. The order we have introduced on $v_G G = v_G(G^{<0})$ is just the one induced by the order of $G^{<0}$.

The maps $\ell$, $\chi_w$ and $\zeta_w$ are monotone. Through the above definition, they induce corresponding equivalence relations on $K^{>0} \setminus R_v$, $(wK)^{<0}$ and $v_G wK$. Since we assume $\ell$ to be strong, the orientation of these maps is fixed (cf. (5), (14) and (15)). Therefore, we have: for $a, a' \in K^{>0} \setminus R_v$, $a \sim_\ell a'$ if and only if there is some $n \in \mathbb{N}$ such that $\ell^n a \leq a'$ and $\ell^n a' \leq a$. Similarly, if $f$ is a strong exponential, then $a \sim_f a'$ if and only if there is some $n \in \mathbb{N}$ such that $f^n a \geq a'$ and $f^n a' \geq a$. If $\ell = f^{-1}$, then the relations $\sim_\ell$ and $\sim_f$ coincide since $f^n a \geq a'$ holds if and only if $a \geq \ell^n a'$ holds. If $g, g' \in (wK)^{<0}$, then $g \sim_{\chi_w} g'$ if and only if there is some $n \in \mathbb{N}$ such that $\chi_w^n g \geq g'$ and $\chi_w^n g' \geq g$. Similarly, if $\gamma, \gamma' \in v_G wK$, then $\gamma \sim_{\zeta_w} \gamma'$ if and only if there is some $n \in \mathbb{N}$ such that $\zeta_w^n \gamma \geq \gamma'$ and $\zeta_w^n \gamma' \geq \gamma$.

**Lemma 9** *For all $g, g' \in (wK)^{<0}$, $g \sim_{\chi_w} g'$ holds if and only if $v_G g \sim_{\zeta_w} v_G g'$ holds. For all $a, a' \in K^{>0} \setminus R_v$, the assertions $a \sim_\ell a'$, $va \sim_\chi va'$ and $v_G va \sim_\zeta v_G va'$ are equivalent.*

Proof:  Suppose that $g \sim_{\chi_w} g'$ and take $n \in \mathbb{N}$ such that $\chi_w^n g \geq g'$ and $\chi_w^n g' \geq g$. Then $\zeta_w^n v_G g = v_G \chi_w^n g \geq v_G g'$ and $\zeta_w^n v_G g' = v_G \chi_w^n g' \geq v_G g$, that is, $v_G g \sim_{\zeta_w} v_G g'$. For the converse, suppose the latter and take $n \in \mathbb{N}$ such that $v_G \chi_w^n g = \zeta_w^n v_G g \geq v_G g'$ and $v_G \chi_w^n g' = \zeta_w^n v_G g' \geq v_G g$. By (13) and (14), this implies that $\chi_w^{n+1} g \geq \chi_w g' \geq g'$ and $\chi_w^{n+1} g' \geq \chi_w g \geq g$, that is, $g \sim_{\chi_w} g'$.

For the proof of our second assertion, it remains to show that $a \sim_\ell a'$ and $va \sim_\chi va'$ are equivalent. Suppose that $a \sim_\ell a'$ and take $n \in \mathbb{N}$ such that $\ell^n a \leq a'$ and $\ell^n a' \leq a$. Then $\chi^n wa = w\ell^n a \geq wa'$ and $\chi^n wa' = w\ell^n a' \geq wa$, that is, $wa \sim_{\chi_w} wa'$. In particular, we obtain that $va \sim_\chi va'$. For the converse, assume the latter and take $n \in \mathbb{N}$ such that $\chi^n va \geq va'$ and $\chi^n va' \geq va$. By (14), we obtain that $v\ell^{n+1} a = \chi^{n+1} va > \chi^n va \geq va'$ and $v\ell^{n+1} a' = \chi^{n+1} va' > \chi^n va' \geq va$. Consequently, $\ell^{n+1} a < a'$ and $\ell^{n+1} a' < a$, that is, $a \sim_\ell a'$.  □

**Lemma 10** *$\sim_\ell$ is coarser than the archimedean equivalence relations on $K^{>0} \setminus R_v$ with respect to addition and multiplication, and $\sim_\chi$ is coarser than the archimedean equivalence* *relation on $(vK)^{<0}$. In other words, the equivalence classes of $\sim_\ell$ are closed under addition and multiplication, and those of $\sim_\chi$ are closed under addition.*

Proof: Assume that $a, a' \in K^{>0} \setminus R_v$ such that $a < a' < na$. Since $\ell$ is strong, we have that $va = vna < v\ell na$ and thus, $\ell a' < \ell na < a$. This proves that archimedean equivalence with respect to addition implies $\ell$-equivalence. Now if $a < a' < a^n$, then $\ell a < \ell a' < n\ell a$, and by what we have already shown, $\ell a \sim_\ell \ell a'$. Since $\ell b \sim_\ell b$ for every $b \in K^{>0} \setminus R_v$, it follows that $a \sim_\ell \ell a \sim_\ell \ell a' \sim_\ell a'$. This proves that archimedean equivalence with respect to multiplication implies $\ell$-equivalence. In view of Lemma 9 and the fact that $v(ab) = va + vb$, this result also yields our assertion about $\sim_\chi$. □

By Lemma 9, $v$ induces an order reversing bijection between $(K^{>0} \setminus R_v)/\sim_\ell$ and $(vK)^{<0}/\sim_\chi$, and $v_G$ induces an order preserving bijection between $(wK)^{<0}/\sim_{\chi_w}$ and $v_G wK/\sim_{\zeta_w}$. In passing, we note (cf. [A], [KS]):

**Theorem 11** *Assume that $\ell$ is surjective. Then $vK$ is divisible and $v_G vK$ is dense without endpoints. With the induced ordering, also every $\chi$-equivalence class and every $\zeta$-equivalence class is dense without endpoints.*

Proof: If $\ell$ is surjective, then so are $\chi$ and $\zeta$. In view of (14) and (15), this yields that their equivalence classes have no endpoints. Further, $\ell^{-1}$ is an exponential and it follows that the positive multiplicative group of $K$ is divisible like its additive group. Consequently, also $vK$ is divisible, hence dense without endpoints (since it is nontrivial by our general assumption that $K$ is nonarchimedean ordered). To show that $v_G vK$ is dense, let $\gamma, \gamma' \in v_G vK$ such that $\gamma < \gamma'$. Using the surjecitivity of $\chi$, we choose $g, g' \in (vK)^{<0}$ with $v_G \chi g = \gamma$ and $v_G \chi g' = \gamma'$. Then $\chi g < \chi g'$, so by density of $vK$ we can choose some $g''$ such that $\chi g < \chi g'' < \chi g'$. In view of (13), this yields that $\gamma = v_G g < v_G g'' < v_G g' = \gamma'$, showing that also $v_G vK$ is dense. Since the $\chi$-equivalence classes and the $\zeta$-equivalence classes are convex subsets of $vK$ and $v_G vK$, it follows that they are also dense. □

We set

$$\Gamma := v_G vK = v_G((vK)^{<0}) \quad \text{and} \quad \Gamma_w := v_G G_w = v_G(G_w^{\leq 0})$$

for every convex subgroup $G_w$ of $vK$. Since $G_w^{\leq 0}$ is a final segment of $(vK)^{<0}$ and $v_G$ preserves $\leq$, we find that $\Gamma_w$ is a final segment of $\Gamma$. We note that

$$a \in R_w \Leftrightarrow va \in G_w \Leftrightarrow v_G va \in \Gamma_w \quad \text{for all } a \in K^{>0} \setminus R_v \qquad (16)$$

(the second equivalence and the implication $va \in G_w \Rightarrow a \in R_w$ hold more generally for all $a \in K$). Indeed, the implication $va \in G_w \Rightarrow v_G va \in \Gamma_w$ holds by definition of $\Gamma_w$. The converse holds since the convex subgroup $G_w$ is closed under archimedean equivalence. Further, we note that $va = vb$ implies that $wa = wb$. Hence, $va \in G_w$ implies that $wa = 0$, whence $a \in R_w$. For the converse, observe that every $a \in K^{>0} \setminus R_v$ satisfies $va < 0$ and thus, $wa \leq 0$. If in addition $a \in R_w$, then $a \in \mathcal{U}_w^{>0}$ and consequently, $va \in G_w$.

**Theorem 12** *Let $\ell$ be a strong logarithm, compatible with $v$. Then the following assertions are equivalent:*

a) $\ell(\mathcal{U}_w^{>0}) = R_w \cap \operatorname{im}(\ell)$      b) $\ell(K^{>0} \setminus R_w) \subset K^{>0} \setminus R_w$

c) $(vK)^{<0} \setminus G_w$ *is closed under* $\chi$      d) $G_w^{<0}$ *is closed under $\chi$-equivalence*

e) $\Gamma \setminus \Gamma_w$ *is closed under* $\zeta$      f) $\Gamma_w$ *is closed under $\zeta$-equivalence.*

*If $\ell$ is the inverse of an exponential $f$, then these conditions are also equivalent to*

g) $a \in R_w \Rightarrow f(a) \in R_w$ *for all* $a \in K^{>0} \setminus R_v$

h) $va \in G_w \Rightarrow vf(a) \in G_w$ *for all* $a \in K^{>0} \setminus R_v$

i) $v_G va \in \Gamma_w \Rightarrow v_G vf(a) \in \Gamma_w$ *for all* $a \in K^{>0} \setminus R_v$.

Proof:    a) $\Leftrightarrow$ b): This was already shown in the last section.

b) $\Leftrightarrow$ c): We know from the last section that condition b) is equivalent to (11). But $wa < 0$ is equivalent to $va < G_w$, and $w\ell a < 0$ is equivalent to $\chi va = v\ell a < G_w$. Thus, (11) is equivalent to condition c).

c) $\Rightarrow$ d): Suppose that $g \in G_w^{<0}$ and that $g' \in (vK)^{<0}$ with $g \sim_\chi g'$. Take $n \in \mathbb{N}$ such that $\chi^n g' \geq g$. Since $G_w^{<0}$ is a final segment of $(vK)^{<0}$, it follows that $\chi^n g' \in G_w^{<0}$. Since $(vK)^{<0} \setminus G_w$ is assumed to be closed under $\chi$, this implies that $g' \in G_w^{<0}$.

d) $\Rightarrow$ c): Take $g \in (vK)^{<0} \setminus G_w$. Since $\chi g \sim_\chi g$ and $G_w^{<0}$ is assumed to be closed under $\chi$-equivalence, we find that $\chi g \in (vK)^{<0} \setminus G_w$.

c) $\Leftrightarrow$ e): Follows directly from the definition of $\zeta$.

e) $\Leftrightarrow$ f): Similar to the proof of c) $\Leftrightarrow$ d).

Now suppose that $\ell = f^{-1}$.

d) $\Rightarrow$ h): Follows from the fact that $va = v\ell f(a) = \chi vf(a)$ is $\chi$-equivalent to $vf(a)$.

h) $\Rightarrow$ d): Suppose that $g \in G_w^{<0}$ and that $g' \in (vK)^{<0}$ with $g \sim_\chi g'$. Choose $a' \in K^{>0}$ such that $g' = va'$. Take $n \in \mathbb{N}$ such that $\chi^n g' \geq g$. Since $G_w^{<0}$ is a final segment of $(vK)^{<0}$, it follows that $v\ell^n a' = \chi^n g' \in G_w^{<0}$. By $n$-fold application of g), we find that $g' = va' = vf^n \ell^n a' \in G_w$.

g) $\Leftrightarrow$ h) $\Leftrightarrow$ i): Follows from (16).      □

Now we are able to give the

**Proof of Theorem 1:** Let $f$ be a strong $T_1$-exponential. Suppose first that $w$ is compatible with $f$. By our remark following (16) and by the first condition of (CE),

$$va \in G_w \Rightarrow a \in R_w \Rightarrow f(a) \in \mathcal{U}_w^{>0} \Rightarrow wf(a) = 0 \Rightarrow vf(a) \in G_w.$$

Now suppose that (1) holds. Then in particular, condition h) of the foregoing theorem holds, which proves that $w$ is compatible with $f$.      □

## 4   Exponential rank and principal exponential rank

Let $(M, <)$ be any totally ordered set. Then the set $M^{\text{fs}}$ of nonempty final segments of $M$ is ordered by inclusion, and the map

$$\iota : M \ni s \mapsto \{s' \in M \mid s' \geq s\}$$

is an order reversing embedding. Its image consists of all segments which have a smallest element.

Assume in addition that $\sim$ is an equivalence relation on $M$ with convex equivalence classes. Let $S, S'$ be final segments of $M$. Then the closure $\tilde{S}$ of $S$ under $\sim$ is again a final segment of $M$. We write $S \sim S'$ if $\tilde{S} = \tilde{S}'$. This definition is compatible with $\iota$ since $s \sim s' \Rightarrow \widetilde{\iota s} = \widetilde{\iota s'} \Rightarrow \iota s \sim \iota s'$. We denote the equivalence classes of $s$ and $S$ by $[s]_\sim$ and $[S]_\sim$. If $M/\sim$ is endowed with the order induced by $<$ and $M^{\text{fs}}/\sim$ is endowed with the order induced by the order $\subset$ of $M^{\text{fs}}$, then the map

$$\sigma: \ M^{\text{fs}}/\sim \ \ni \ [S]_\sim \ \mapsto \ \{[s]_\sim \mid s \in S' \in [S]_\sim\} \ = \ \{[s]_\sim \mid s \in S\} \ \in \ (M/\sim)^{\text{fs}} \qquad (17)$$

is an order preserving bijection.

Now let $g, g' \in vK$. Denote by $C_g$ the smallest convex subgroup of $vK$ containing $g$. Then $C_g \subsetneq C_{g'}$ if and only if $|g'| > n|g|$ for all $n \in \mathbb{N}$, i.e., if and only if $v_G g' < v_G g$. Hence, the map

$$\{C_g \mid 0 \neq g \in vK\} \ni C_g \mapsto v_G g \in \Gamma$$

is an order reversing bijection. Composing this bijection with the map $R_w \mapsto G_w$, we obtain an order reversing bijection

$$\mathcal{R}^{\text{pr}} \to \Gamma \ . \qquad (18)$$

Every convex subgroup $G_w$ is the union of the principal convex subgroups contained in it. Correspondingly, every $R_w \in \mathcal{R}$ is the union of all rings in $\mathcal{R}^{\text{pr}}$ which are contained in $R_w$. Such a union corresponds via (18) to the final segment $\Gamma_w = v_G G_w$ of $\Gamma = v_G vK$. We have thus obtained an order preserving bijection

$$\rho: \ \mathcal{R} \ni R_w \mapsto \Gamma_w \in \Gamma^{\text{fs}} \ .$$

(As we have excluded $R_v$ from $\mathcal{R}$, we do not have to deal with $v_G 0 = \infty$ which by our definition does not lie in $v_G vK$.) Let us identify the set $\Gamma$ with its image in $\Gamma^{\text{fs}}$ (which consists of all segments having a smallest element). Then the restriction of the map $\rho$ to $\mathcal{R}^{\text{pr}}$ is just the bijection (18). Since $\zeta$ preserves $\leq$ on $\Gamma$, it sends final segments to final segments. That is, $\zeta$ and thus also the equivalence relation $\sim_\zeta$ extend canonically to $\Gamma^{\text{fs}}$. Via the bijection $\rho$, we may view $\zeta$ as a map on $\mathcal{R}$.

Now we consider the following map, where $\sigma$ is defined as in (17):

$$\varepsilon: \ \mathcal{R} \ni R_w \mapsto \sigma([\Gamma_w]_{\sim_\zeta}) \in (\Gamma/\sim_\zeta)^{\text{fs}} \ .$$

It is an epimorphism which preserves $\subseteq$. To obtain a bijection, we restrict our scope to $R_w \in \mathcal{R}_\ell$: since the corresponding final segments $\Gamma_w$ of $\Gamma$ are closed under $\zeta$-equivalence, the map $\Gamma_w \mapsto [\Gamma_w]_{\sim_\zeta}$ becomes injective. Hence, we obtain an order preserving bijection $\varepsilon$ from $\mathcal{R}_\ell$ onto $(\Gamma/\sim_\zeta)^{\text{fs}}$. We compute:

$$\begin{aligned}
\varepsilon(R_w) &= \{[\gamma]_\zeta \mid \gamma \in \Gamma_w\} = \{[v_G g]_\zeta \mid 0 \neq g \in G_w\} \\
&= \{[v_G va]_\zeta \mid a \in K^{>0} \wedge wa = 0 \wedge va \neq 0\} \\
&= \{[v_G va]_\zeta \mid a \in \mathcal{U}_w^{>0} \wedge va \neq 0\} = \{[v_G va]_\zeta \mid a \in \mathcal{U}_w^{>0} \setminus \mathcal{U}_v^{>0}\} \ .
\end{aligned}$$

We summarize what we have proved:

**Theorem 13** *Let $\ell$ be a (not necessarily surjective) strong logarithm. Then*

$$\varepsilon: R_w \mapsto \{[v_G va]_\zeta \mid a \in \mathcal{U}_w^{>0} \setminus \mathcal{U}_v^{>0}\}$$

*is an order preserving bijection from the exponential rank onto $(\Gamma/\sim_\zeta)^{\text{fs}}$.*

We define the **principal exponential rank** $\mathcal{R}_\ell^{\text{pr}}$ to be the preimage under $\varepsilon$ of the set of all final segments of $\Gamma/\sim_\zeta$ which have a smallest element. Observe that $\varepsilon(R_w)$ contains a smallest element if and only if $\Gamma_w$ admits some $\sim_\zeta$-equivalence class as initial segment, or equivalently, $G_w^{<0}$ admits some $\chi$-equivalence class as initial segment. This does not mean that $G_w$ is principal; the following corollary shows the contrary.

**Corollary 14** *If $\ell$ is surjective, then the intersection of the principal rank and the exponential rank is empty. In particular, the value group of a nonarchimedean exponential field is never principal (as its own convex subgroup).*

Proof: If $R_w$ belongs to the exponential rank, then $v_G G_w$ is closed under $\zeta$-equivalence. If $\ell$ is surjective, then Theorem 11 shows that $\zeta$-equivalence classes have no smallest element; hence also $v_G G_w$ has no smallest element. The second assertion follows from the first, taking $w$ to be the trivial valuation. □

For the rest of this chapter, assume that $\ell$ is the inverse of a strong $T_1$-exponential $f$. First, we wish to show that $\mathcal{R}_\ell^{\text{pr}} = \mathcal{R}_f^{\text{pr}}$. By Theorem 12 and Theorem 1, $G_w$ is closed under $va \mapsto vf(a)$ if and only if $\Gamma_w$ is closed under $\zeta$-equivalence. Therefore, $G_w$ is the smallest convex subgroup containing $g$ and closed under $va \mapsto vf(a)$ if and only if $\Gamma_w$ is the smallest final segment containing $v_G g$ and closed under $\zeta$-equivalence. This in turn holds if and only if $[v_G g]_\zeta$ is an initial segment of $\Gamma_w$. Consequently, $\mathcal{R}_\ell^{\text{pr}} = \mathcal{R}_f^{\text{pr}}$.

Take any $a \in K^{>0} \setminus R_v$ and $R_w \in \mathcal{R}_f$, $w \neq v$. From the results of the previous chapter it follows that the sequence $(f^n a)_{n \in \mathbb{N}}$ is cofinal in the $\sim_f$-equivalence class $[a]_f$ of $a$, the sequence $(vf^n a)_{n \in \mathbb{N}}$ is coinitial in the $\sim_\chi$-equivalence class $[va]_\chi$ of $va$, and the sequence $(v_G v f^n a)_{n \in \mathbb{N}}$ is coinitial in the $\sim_\zeta$-equivalence class $[va]_\zeta$ of $v_G va$. Hence, the sequence $(f^n a)_{n \in \mathbb{N}}$ is cofinal in $R_w$ if and only if $[va]_\chi$ is an initial segment of $G_w$, and this is the case if and only if $[v_G va]_\zeta$ is an initial segment of $\Gamma_w$. This in turn holds if and only if $[v_G va]_\zeta$ is the smallest element in $\varepsilon(R_w)$.

On the other hand, the sequence $(f^n a)_{n \in \mathbb{N}}$ is cofinal in $R_w$ if and only if the sequence $((fw)^n(aw))_{n \in \mathbb{N}}$ is cofinal in $Kw$. This follows from $(f^n a)w = (fw)^n(aw)$ and the fact that the residue map $a \mapsto aw$ induces a $\leq$-preserving map from $R_w^{>0} \setminus R_v^{>0}$ onto the positive infinite elements of $Kw$. We have thus proved the first assertion of part b) of Theorem 2.

If for every $a \in K^{>0} \setminus R_v$ the sequence $(f^n a)_{n \in \mathbb{N}}$ is cofinal, then this means that for every such $a$ the class $[v_G va]_\zeta$ is the same, and vice versa. This in turn means that $\Gamma/\sim_\zeta$ is a singleton, i.e., the principal exponential rank is a singleton. This proves the second assertion of part b) of Theorem 2. The proof of part a) is similar.

Finally, let us mention (and leave the proof as an exercise to the reader):

**Theorem 15** *Assume that $w$ is compatible with $f$. Then $\chi$ induces $\chi_w$ through the canonical isomorphism $wK \simeq vK/G_w$, and $\zeta_w$ is the restriction of $\zeta$ to $v_G wK$ through the* *canonical isomorphism* $v_G wK \simeq v_G vK \setminus \Gamma_w$. *Further, the valuation $\overline{w}$ induced by $v$ on the residue field $Kw$ is the natural valuation of $Kw$ (endowed with the induced order), the exponential $fw$ on $Kw$ induces the restriction of $\chi$ on the value group $G_w$ through the canonical isomorphism $\overline{w}(Kw) \simeq G_w$, and the restriction of $\zeta$ on $\Gamma_w$ through the canonical isomorphism $v_G \overline{w}(Kw) \simeq \Gamma_w$.*

## 5 Lifting logarithms from the residue field

If $K$ admits an exponential, then its multiplicative group of positive elements is divisible (since the additive is). For the rest of the paper, we will always assume this divisibility. As in [KS] (Lemma 3.4 and Theorem 3.8), we then have the following representations as lexicographic sums:

$$(K, +, 0, <) \simeq \mathbf{A}_w \amalg (Kw, +, 0, <) \amalg (I_w, +, 0, <) \tag{19}$$

where $\mathbf{A}_w$ is an arbitrary group complement of $R_w$ in $(K, +)$, and analogously,

$$(K^{>0}, \cdot, 1, <) \simeq \mathbf{B}_w \amalg (Kw^{>0}, \cdot, 1, <) \amalg (1 + I_w, \cdot, 1, <) \tag{20}$$

where $\mathbf{B}_w$ is an arbitrary group complement of $\mathcal{U}_w^{>0}$ in $(K^{>0}, \cdot)$. Endowed with the restriction of the ordering, $\mathbf{A}_w$ and $\mathbf{B}_w$ are unique up to isomorphism. In view of (CO) and the fact that $w(-a) = wa$, the map

$$(K^{>0}, \cdot, 1, <) \to (wK, +, 0, <), \qquad a \mapsto -wa = wa^{-1} \tag{21}$$

is a surjective group homomorphism preserving $\leq$, with kernel $\mathcal{U}_w^{>0}$. We find that every complement $\mathbf{B}_w$ is isomorphic to $(wK, +, 0, <)$ through the map $-w$.

Let $w$ be compatible with the (not necessarily surjective) logarithm $\ell$. Then $\ell$ decomposes into three embeddings of ordered groups:

$$\begin{array}{rrcl} \ell_R^w : & (1 + I_w, \cdot, 1, <) & \to & (I_w, +, 0, <) \\ \ell w : & (Kw^{>0}, \cdot, 1, <) & \to & (Kw, +, 0, <) \\ \ell_L^w : & \mathbf{B}_w & \to & \mathbf{A}_w \,. \end{array}$$

Conversely, in view of (19) and (20), such embeddings $\ell_R^w$, $\ell w$ and $\ell_L^w$ can be put together to obtain a logarithm which is compatible with $w$. We call $\ell_L^w$ a **left logarithm** and $\ell_R^w$ a **right logarithm**. $\ell w$ is a logarithm on the residue field $Kw$, and $\ell$ can be seen as a **lifting** of $\ell w$. Thus, the liftings of $\ell w$ to $K$ are in one-one correspondance to the pairs $(\ell_L^w, \ell_R^w)$ of left and right logarithms. The set of all right logarithms is identical to the set of all order preserving embeddings of $(1 + I_w, \cdot)$ in $(I_w, +)$; we will denote it by o-Emb $((1 + I_w, \cdot), (I_w, +))$.

Through the isomorphism (21), every embedding

$$h : (wK, +, 0, <) \to \mathbf{A}_w$$

gives rise to a left logarithm $h \circ -w$. Conversely, given a left logarithm $\ell_L^w$, the map

$$h_\ell^w := \ell_L^w \circ (-w)^{-1}$$

is such an embedding (here, $(-w)^{-1}$ is an isomorphism from $wK$ onto $\mathbf{B}_w$); note that $h$ is surjective if and only if $\ell_L^w$ is. This one-to-one correspondence motivates the following definition. A **logarithmic cross-section** of an ordered field $(K,<)$ with respect to a convex valuation $w$ is an order preserving embedding $h$ of $wK$ into an additive group complement of the valuation ring, or equivalently, an embedding $h$ of $wK$ in $(K,+,0,<)$ such that $h(wK) \cap R_w = \{0\}$; we will denote the set of all such embeddings by o-Emb$(wK, (K,+) \setminus R_w)$.

Further, we denote the set of all (not necessarily surjective) logarithms of $K$ by $\mathbf{L}_K$, and $\mathbf{L}_K^w$ shall be the subset of those logarithms which are compatible with $w$. Then we have:

**Theorem 16** *The map*

$$\begin{array}{rcl}
\mathbf{L}_K^w & \to & \text{o-Emb}(wK, (K,+) \setminus R_w) \times \mathbf{L}_{Kw} \times \text{o-Emb}((1+I_w, \cdot), (I_w, +)) \\
\ell & \mapsto & (h_\ell^w, \ell w, \ell_R^w)
\end{array} \quad (22)$$

*is a bijection, and the following holds:*
*a) $\ell$ is surjective if and only if $h_\ell^w$, $\ell w$ and $\ell_R^w$ are,*
*b) if $w'$ is a coarsening of $v$ such that $w$ is a coarsening of $w'$, then $\ell$ is compatible with $w'$ if and only if $\ell w$ is compatible with the induced valuation $w'/w$.*

(Use a lexicographic decomposition of $Kw$ similar to the above to prove the last assertion.)

Let us quickly compare this result with the lifting of orderings through places. If we denote by $\mathbf{X}_K$ the set of all orderings on $K$, and $\mathbf{X}_K^w$ the subset of all orderings which are compatible with $w$, then there is a bijection

$$\mathbf{X}_K^w \to \text{Hom}(vK/2vK, \{-1,1\}) \times \mathbf{X}_{Kw}.$$

We wish to derive a condition for $\ell$ to be strong. Every $a \in K^{>0} \setminus R_w$ can be written as $a = b \cdot c$ where $b \in \mathbf{B}_w$ and $c \in \mathcal{U}_w^{>0}$, and $wa = wb$. Then $w\ell a = w(\ell b + \ell c) = w\ell b$ since $b \in \mathbf{B}_w$, $c \in \mathcal{U}_w^{>0}$ imply that $w\ell b < 0 \leq w\ell c$. Hence,

$$wa < w\ell a \qquad \text{for all } a \in K^{>0} \setminus R_w$$

is equivalent to

$$w\ell_L^w a > wa \qquad \text{for all } a \in \mathbf{B}_w^{>0}. \quad (23)$$

With $g = wa$ and the isomorphism $h = \ell_L^w \circ (-w)^{-1} : wK \to \mathbf{A}_w$, we have that $w\ell_L^w a = w(-\ell_L^w a) = w\ell_L^w(a^{-1}) = w\ell_L^w \circ (-w)^{-1} \circ (-w)(a^{-1}) = wh(wa)$. Hence, condition (23) translates to

$$wh(g) > g \qquad \text{for all } g \in (wK)^{<0}. \quad (24)$$

If $h \in$ o-Emb$(wK, (K,+) \setminus R_w)$ satisfies $wh(g) > g$ for all $g \in (wK)^{<0}$, then we call it a **strong logarithmic cross-section** (for $w$). For $w = v$, we see that (6) holds if and only if (24) holds for $h_\ell^v$. We have thus proved the first part of the following lemma, and we leave the proof of the second part to the reader:

**Lemma 17** *A logarithm $\ell$ is strong if and only if $h_\ell^v$ is strong. If $\ell \in \mathbf{L}_K^w$ and $\ell w$ and $h_\ell^w$ are strong, then also $\ell$ is strong.*

The converse of the last assertion does not hold in general. If $h_\ell^w$ is strong, then "$g < \chi_w g$" holds in (14) and "$\gamma < \zeta_w \gamma$" holds in (15). If $\ell$ is surjective and $w \neq v$, this describes a more rapid growth rate of the exponential $\ell^{-1}$ on the positive infinite elements than the axiom (GA) does.

Theorem 16 does not yet tell anything about the existence of (strong) logarithmic cross-sections and right logarithms (if we don't know whether logarithms exist). We will now discuss this problem. Recall that every embedding (resp. isomorphism) of ordered abelian groups induces canonically an embedding (resp. isomorphism) of their ranks as ordered sets (cf. [KS1]). In particular, a logarithmic cross-section $h$ induces an embedding $\tilde{h}$ such that the following diagram commutes:

$$\begin{array}{ccc} vK & \xrightarrow{h} & \mathbf{A}_v \\ {\scriptstyle v_G} \downarrow & & \downarrow {\scriptstyle v} \\ v_G vK & \xrightarrow{\tilde{h}} & (vK)^{<0} \end{array}$$

We say that $h$ is a **lifting** of $\tilde{h}$. If $h$ is onto, then so is $\tilde{h}$ (in this case, it is just the inverse of a "group exponential" as defined in [KS1]). We have that

$$\tilde{h}(v_G g) > g \Leftrightarrow vh(g) > g$$

for every $g \in (vK)^{<0}$.

We see that $h$ is a strong logarithmic cross-section if and only if

$$\tilde{h}(v_G g) > g \qquad \text{for all } g \in (vK)^{<0} . \tag{25}$$

Note that every ordered abelian group $G$ admits an embedding $s : v_G G \to G^{<0}$ of ordered sets such that $v_G \circ s$ is the identity on $v_G G$ (for $\alpha \in v_G G$, we just have to set $s\alpha = g$ where $g \in G^{<0}$ is an arbitrary element of value $v_G g = \alpha$). We call such a map a **group cross-section**.

**Lemma 18** *Let $G$ be any ordered abelian group such that $v_G G$ admits an order preserving map $\zeta$ into itself satisfying that $\zeta \alpha > \alpha$ for all $\alpha \in v_G G$. Then for every group cross-section $s$ of $G$, the embedding $\tilde{h} := s \circ \zeta : v_G G \to G^{<0}$ will satisfy condition (25).*

Indeed, $v_G \tilde{h}(v_G g) = \zeta v_G g > v_G g$ and thus a fortiori $\tilde{h}(v_G g) > g$ if $g \in G^{<0}$. Note that there are plenty of groups satisfying the hypothesis of the lemma. For instance, this is the case if $v_G G$ is isomorphic to an arbitrary nontrivial ordered abelian group, as an ordered set.

Now the question arises whether an embedding (resp. isomorphism) $\tilde{h}$ can be lifted to an embedding (resp. isomorphism) $h$. (Cf. the related notion of "lifting property" as used in [KK1].) Such a lifting always exists if $\mathbf{A}_v$ is rich enough, i.e., if it is a Hahn product. This in turn is the case if the field $K$ is a suitable power series field.

Let $k$ be an archimedean ordered real closed field. If $G$ is an arbitrary ordered abelian group, then the power series field $K := k((G))$ is a formally real field, and it is real closed if and only if $G$ is divisible (which we shall always assume here). Further, $K$ carries a canonical valuation $v$ which associates to every formal power series the minimum of its support. It also carries a natural ordering $<$ such that $v$ is the natural valuation of the ordered field $(K,<)$. The residue field of $(K,v)$ is $k$, and its value group is $G$. The valuation ring $R$ of $(K,v)$ is the power series ring $k[[G]]$. We can take the additive group complement $\mathbf{A}_v$ of the valuation ring $R_v$ to be the ordered ring $k((G^{<0})) := \{a \in k((G)) \mid \text{support}(a) \subset G^{<0}\}$. As an ordered abelian group, it is canonically isomorphic to the Hahn product $\mathbf{H}_{G^{<0}}(k,+,0,<)$. Concerning the existence of right exponentials, the following result is well known:

**Lemma 19** *Let $I_v$ be the valuation ideal of $k((G))$. Then for every $\varepsilon \in I_v$,*

$$f_R(\varepsilon) := \sum_{i=0}^{\infty} \frac{\varepsilon^i}{i!} \quad (26)$$

*is a canonically defined element of $1 + I_v$ (cf. Neumann's Lemma [N]), and $\ell_R := f_R^{-1}$ is a surjective right logarithm. "Canonical" means in particular: if $G \subset G'$ and $f'_R$ is defined on the valuation ideal of $k((G'))$ in the same way, then it extends $f_R$.*

For the case of $k = \mathbb{R}$, we can show:

**Theorem 20** *Let $\tau$ be any order type. Then there is a divisible ordered abelian group $G$ such that $v_G G$ admits an automorphism $\zeta$ satisfying*
*a) $\zeta\alpha > \alpha$ for all $\alpha \in v_G G$*
*b) $v_G G/\sim_\zeta$ has order type $\tau$.*
*Further, the power series field $\mathbb{R}((G))$ admits a strong logarithmic cross-section for $v$, giving rise to a (non-surjective) strong logarithm having principal exponential rank $\tau$.*

Proof: Let $T$ be an ordered set having order type $\tau$. We may assume that $\tau$ is nontrivial, that is, $T \neq \emptyset$, since otherwise, we could set $G = \{0\}$ and $\mathbb{R}((G)) = \mathbb{R}$, and the usual logarithm would do the job. We define the ordered set $\Gamma$ to be the sum (in the sense of ordered sets) of copies of $\mathbb{Z}$ over the index set $T$. (That is, we obtain $\Gamma$ by replacing every element of $T$ by a copy of $\mathbb{Z}$). We let $\zeta$ be the map which sends an element $n$ in any of these copies to its successor $n+1$ in the same copy. Now we let $G$ be the Hahn sum (or Hahn product) of copies of $\mathbb{Q}$ over the index set $\Gamma$. Then $G$ has the required properties.

According to Lemma 18, we can choose an embedding $\tilde{h} : v_G G \to G^{<0}$ which satisfies condition (25). Note that $\mathbf{A}_v$ is archimedean-complete (that is, it is maximal and all its components are $\mathbb{R}$). Hence by Hahn's embedding theorem, the embedding $\tilde{h}$ of $v_G G$ into $G^{<0} = v(\mathbf{A}_v \setminus \{0\})$ lifts to an embedding $h$ of $G$ into $\mathbf{A}_v$. Moreover, since $\tilde{h}(v_G g) > g$, we have that $vh(g) > g$ for all $g \in G^{<0}$, as required.

In view of the foregoing lemma, Theorem 16 now shows that $h$ gives rise to a strong logarithm $\ell$ which lifts $\exp$ from $\mathbb{R}$ to $\mathbb{R}((G))$. To show that $\ell$ has principal exponential rank $\tau$, it suffices to prove that $\ell$ induces $\zeta$ on $v_G G$. As $\zeta$ is already induced by $\ell_L^v$, we take $a \in \mathbf{B}_v$ and compute:

$$\begin{aligned} v_G v \ell_L^v a &= v_G v(h \circ (-v)(a)) = v_G(v(-h(va))) = v_G(v(h(va))) \\ &= v_G \circ v \circ h(va) = v_G \circ \tilde{h} \circ v_G(va) = v_G \circ s \circ \zeta(v_G va) = \zeta(v_G va) , \end{aligned}$$

as required. This completes our proof. $\square$

The so obtained logarithm can never be surjective. If it were, it would give rise to an exponential on the power series field $\mathbb{R}((G))$, compatible with the natural valuation; but this is impossible by the main result of [KKS].

# 6 Going to the limit

Using Theorem 20, we shall now construct nonarchimedean models of real exponentiation which are countable unions of power series fields. Indeed, a common method to obtain surjectivity of a map is to construct the union over a suitable countably infinite chain of fields. In the following, we will apply such a construction to strong logarithmic cross-sections.

- **Construction of a surjective logarithmic cross-section.**

To get started, let $G$ be as in Theorem 20. Set $G_0 := G$ and $K_0 = \mathbb{R}((G_0))$. Let $\mathbf{A}_0$ be a group complement of $\mathbb{R}[[G_0]]$ in $K_0$ and $h_0 : G_0 \to \mathbf{A}_0$ a strong logarithmic cross-section of $K_0$. Now assume that we have already constructed $G_{n-1}$, $K_{n-1}$, $\mathbf{A}_{n-1}$ and the strong logarithmic cross-section
$$h_{n-1} : G_{n-1} \to \mathbf{A}_{n-1}$$
satisfying
$$vh_{n-1}(g) > g \quad \text{for all } g \in G_{n-1}^{<0} . \tag{27}$$

Since $G_{n-1}$ is isomorphic to a subgroup of $\mathbf{A}_{n-1}$ through $h_{n-1}$, we can take $G_n$ to be a group containing $G_{n-1}$ as a subgroup and admitting an isomorphism $h_n$ onto $\mathbf{A}_{n-1}$ which extends $h_{n-1}$. We set $K_n := \mathbb{R}((G_n))$. Hence, $K_{n-1} \subset K_n$ canonically (the elements of $K_{n-1}$ being those elements of $K_n$ whose support is a subset of $G_{n-1}$). Further, we choose a group complement $\mathbf{A}_n$ for the valuation ring $\mathbb{R}[[G_n]]$ such that $\mathbf{A}_n$ contains $\mathbf{A}_{n-1}$. In this way, $h_n$ appears as an embedding of $G_n$ into $\mathbf{A}_n$ which extends $h_{n-1}$. We show that $h_n$ is again a strong logarithmic cross-section. For $g \in G_n$, the image $h_n(g)$ lies in $\mathbf{A}_{n-1}$, and $vh_n(g)$ lies in its value set $G_{n-1}^{<0}$. Consequently, in (27) we may replace $g \in G_{n-1}^{<0}$ by $vh_n(g)$ for $g \in G_n^{<0}$. But $vh_{n-1}(vh_n(g)) > vh_n(g)$ implies that $h_{n-1}(vh_n(g)) > h_n(g)$, because $h_n(g) < 0$ and $h_{n-1}(vh_n(g)) < 0$. Since $h_n$ extends $h_{n-1}$, this may be read as $h_n(vh_n(g)) > h_n(g)$. Since $h_n$ is order preserving, this in turn implies $vh_n(g) > g$. Thus, we have proved that (27) holds with $n$ in the place of $n-1$.

By our induction on $n$, we obtain a chain of fields $K_n$, $n \in \mathbb{N}$. Now we take $K_\omega := \bigcup_{n \in \mathbb{N}} K_n$ and $h_\omega := \bigcup_{n \in \mathbb{N}} h_n$. Also the groups $G_n$ form a chain, and their union $G_\omega := \bigcup_{n \in \mathbb{N}} G_n$ is the value group of $K_\omega$; we have that $K_\omega \subset \mathbb{R}((G_\omega))$. Similarly, the group complements $\mathbf{A}_n$ form a chain, and their union $\mathbf{A}_\omega := \bigcup_{n \in \mathbb{N}} \mathbf{A}_n$ is a group complement for the valuation ring $\bigcup_{n \in \mathbb{N}} \mathbb{R}[[G_n]] = \mathbb{R}[[G_\omega]] \cap K_\omega$ in $K_\omega$. By construction, we have $\mathbf{A}_{n-1} = h_n(G_n)$ for all $n$. Consequently, $h_\omega : G_\omega \to \mathbf{A}_\omega$ is surjective. Moreover, $h_\omega$ is a strong logarithmic cross-section since (27) holds for all $n$.

- **Construction of a surjective strong exponential $f_\omega$ on $K_\omega$.**

By Lemma 19 we obtain a right exponential $f_{R,n}$ and a surjective right logarithm $\ell_{R,n} = f_{R,n}^{-1}$ on every $K_n$, such that $\ell_{R,n+1}$ is an extension of $\ell_{R,n}$ to $K_{n+1}$. Hence, $\ell_{R,\omega} := \bigcup_{n \in \mathbb{N}} \ell_{R,n}$ is a surjective right logarithm on $K_\omega$.

Now we apply Theorem 16 to find a surjective strong logarithm $\ell_\omega$ which lifts the usual exponential function exp from $\mathbb{R}$ to $K_\omega$. Its inverse $f_\omega$ is a strong exponential on $K_\omega$. This completes our construction.

- **Model theoretic properties of** $(K_\omega, f_\omega)$.

In [DMM1] it is shown how to interprete the restricted analytic functions on power series fields via their Taylor expansions. This interpretation is canonical in the same spirit as in Lemma 19, hence it is compatible with the inclusions $K_n \subset K_{n+1}$. Moreover, it makes every $K_n$ into a model of $T_{\text{an}}$. By the model completeness of $T_{\text{an}}$ (cf. [D1]), $K_n \prec K_{n+1}$ for every $n$. Hence, $K_\omega$ is the union over an elementary chain of models $K_n$ of $T_{\text{an}}$ and is thus itself a model of $T_{\text{an}}$.

For every $n \in \mathbb{N}$, $h_n$, log and $\ell_{R,n}$ give rise to a logarithm $\ell_n$ on $K_n$. Note that $\ell_\omega = \bigcup_{n \in \mathbb{N}} \ell_n$. Since every $\ell_{R,n}$ is surjective and exp is surjective on $\mathbb{R}$, the restriction $\ell_n^{\text{fin}}$ of $\ell_n$ to $\mathcal{U}_v^{>0}$ is an isomorphism onto the valuation ring $\mathbb{R}[[G_n]]$. We denote its inverse by $f_n^{\text{fin}}$; it is the restriction of $f_\omega$ to $\mathbb{R}[[G_n]]$.

Let $n \in \mathbb{N}$ and $a$ be an element of the interval $[-1, 1]$ of $K_n$. Since $[-1, 1] \subset \mathbb{R}[[G_n]]$, we can write $a = r + \varepsilon$ with $r \in \mathbb{R}$ and $v\varepsilon > 0$, and we have:

$$f_\omega(a) = f_n^{\text{fin}}(a) = (\ell_n^{\text{fin}})^{-1}(a) = \log^{-1}(r) \cdot \ell_{R,n}^{-1}(\varepsilon)$$
$$= \exp(r) \cdot f_{R,n}(\varepsilon) = \exp(r) \cdot \sum_{i=0}^{\infty} \frac{\varepsilon^i}{i!}.$$

Therefore, $f_\omega$ coincides on $[-1, 1]$ in $K_n$ with the interpretation of the restricted exp (given by its Taylor expansion), for every $n$. Hence, this is also true on the interval $[-1, 1]$ in $K_\omega$. From Lemma 4 we conclude that $(K_\omega, f_\omega)$ is a model of $T_{\text{an}}(\exp)$.

- **The principal exponential rank of** $(K_\omega, f_\omega)$.

We wish to show that $(K_\omega, f_\omega)$ has the same principal exponential rank as $K_0$ with its logarithm induced by $h_0$. Let $a \in \mathbf{A}_\omega^{>0}$; then there is some $n \in \mathbb{N}$ such that $a \in \mathbf{A}_n^{>0}$. By construction, the image of $h_n$ is $\mathbf{A}_{n-1}$. Consequently, $\ell_\omega a \in \mathbf{A}_{n-1}^{>0}$. By induction on $n$, we find that $\ell_\omega^n a \in \mathbf{A}_0$. Since every infinite positive element in $K_\omega$ is archimedean equivalent (and thus $\ell_\omega$-equivalent) to some $a \in \mathbf{A}_\omega$ and $a$ is $\ell_\omega$-equivalent to $\ell_\omega^n a$, this proves that every infinite positive element in $K_\omega$ is $\ell_\omega$-equivalent to some infinite positive element in $K_0$. This proves our assertion.

**Remark 21** The above construction can be iterated in order to obtain unions over chains indexed by an arbitrary limit ordinal $\kappa$. If $\lambda \leq \kappa$ is a limit ordinal and we have constructed $G_\nu$, $K_\nu$, $\mathbf{A}_\nu$ and $h_\nu$ for every $\nu < \lambda$, then we take for $G_\lambda$, $K_\lambda$ and $h_\lambda$ the respective unions in the same manner as before. If $\lambda < \kappa$, then we replace $K_\lambda$ by $\mathbb{R}((G_\lambda))$, which by virtue of the main result of [KKS] must be a proper extension of $\bigcup_{\nu < \lambda} K_\nu$. We choose a group complement $\mathbf{A}_\lambda$ to its valuation ring $\mathbb{R}[[G_\lambda]]$ which contains $\bigcup_{\nu < \lambda} \mathbf{A}_\nu$. Thus, $h_\lambda$ is a non-surjective logarithmic cross-section of $K_\lambda$ with image in $\mathbf{A}_\lambda$. The induction step for successor ordinals works as before.

If $\kappa$ is an uncountable regular cardinal, then the exponential field $(K_\kappa, f_\kappa)$ obtained by this construction is almost a power series field. In fact, it is the restricted power series field $\mathbb{R}((G_\kappa))_\kappa$, which consists of all power series in $\mathbb{R}((G_\kappa))$ whose support has cardinality $< \kappa$. Indeed, since $\kappa$ is assumed to be regular and $G_\kappa = \bigcup_{\nu < \kappa} G_\nu$, every power series with support of cardinality $< \kappa$ is already an element of $\mathbb{R}((G_\nu)) = K_\nu$ for some $\nu < \kappa$. Hence, it lies in $K_\kappa = \bigcup_{\nu<\kappa} K_\nu$.

Now let $(K, f)$ be a model of $T_{\text{an}}(\exp)$, and $\tau$ an order type extending the principal exponential rank $\tau_0$ of $(K, f)$. By abuse of terminology, we assume $\tau_0$ and $\tau$ to be ordered sets of the respective order types. Now for every element in $\tau \setminus \tau_0$ we add a copy of $\mathbb{Z}$ to $\Gamma = v_G v K$, defining $\zeta$ on this copy to send $n$ to its successor $n + 1$. In this way, we obtain an ordered set $\Delta$ with a map $\zeta$ such that $(\Gamma, \zeta)$ embeds in $(\Delta, \zeta)$ and $\Delta/\sim_\zeta \simeq \tau$. We take $G$ to be the Hahn product of copies of $\mathbb{R}$ over the index set $\Delta$. Then $G$ satisfies properties a) and b) of Theorem 20. By Hahn's embedding theorem, the embedding of $\Gamma$ in $\Delta$ lifts to an embedding of $vK$ in $G$.

By [DMM1], $\mathbb{R}((G))$ is a model of the theory $T_{\text{an}}$ of the reals with restricted analytic functions. Moreover, there is a truncation closed embedding of $K$ in $\mathbb{R}((G))$ which respects the restricted analytic functions. Now the left logarithm of $K$ induces canonically a strong logarithmic cross-section $h_0$ on $K_0 = \mathbb{R}((G))$. We continue the construction as above. The so obtained exponential $f_\omega$ on $K_\omega$ extends $f$. By [DMM1], the embedding of $(K, f)$ in $(K_\omega, f_\omega)$ is elementary. Note that in our construction, every embedding $K_{n-1} \subset K_n$ is truncation closed. Hence, the embedding $K \subset K_\omega$ is truncation closed. This proves Theorem 3.

The Fields Institute  
222 College Street  
Toronto, Ontario M5T 3J1, Canada  
email: fkuhlman@fields.utoronto.ca, skuhlman@fields.utoronto.ca